\documentclass[a4paper,12pt]{amsart}

\usepackage{graphicx}
\usepackage{amsmath}

\textwidth=\paperwidth \advance\textwidth by-2in
\makeatletter
\def\temp{12}
\ifx\@mainsize\temp \headheight=6.5pt \fi
\makeatother

\def\Z{\mathbf{Z}}
\def\Q{\mathbf{Q}}
\def\R{\mathbf{R}}

\def\lk{\operatorname{lk}}

\def\inte{\operatorname{int}}

\def\spmatrix#1{\left[
  \begin{smallmatrix} #1 \end{smallmatrix}
  \right]\ignorespaces}
\theoremstyle{plain}
\newtheorem{thm}{Theorem}[section]

\newtheorem{lem}{Lemma}[section]
\theoremstyle{definition}

\theoremstyle{remark}

\begin{document}

\title{Signatures of Covering Links}
\date{August 28, 2001 (First Edition: November 23, 1999)}
\def\subjclassname{\textup{2000} Mathematics Subject Classification}
\expandafter\let\csname subjclassname@1991\endcsname=\subjclassname
\expandafter\let\csname subjclassname@2000\endcsname=\subjclassname
\author{Jae Choon Cha}
\email{jccha\char`\@knot.kaist.ac.kr}
\address{Department of Mathematics\\
Korea Advanced Institute of Science and Technology\\
Taejon, 305--701\\
Korea}
\author{Ki Hyoung Ko}
\email{knot\char`\@knot.kaist.ac.kr}
\address{Department of Mathematics\\
Korea Advanced Institute of Science and Technology\\
Taejon, 305--701\\
Korea}

\begin{abstract}
The theory of signature invariants of links in rational homology spheres is
applied to covering links of homology boundary links. From patterns and
Seifert matrices of homology boundary links, an explicit
formula is derived to compute signature invariants of their covering links.
Using the formula, we produce fused boundary links
that are positive mutants of ribbon
links but are not concordant to boundary links.
We also show that for any finite
collection of patterns, there are homology boundary links that are not
concordant to any homology boundary links admitting a pattern in the
collection.
\end{abstract}

\maketitle

\section{Introduction}

For a link $L$, the pre-image (of a sublink) of $L$ in a finite cyclic
cover of the ambient space branched along a component of $L$ is called
a \emph{covering link of~$L$}.  In the work of Cochran and
Orr~\cite{CO1,CO2}, it was observed that concordances of links in
spheres can be studied via their covering links due to the following
facts: If $L$ is a link in a $\Z_p$-homology sphere for some prime
$p$, so is a $p^a$-fold covering link~\cite{CG2}, and corresponding
covering links of concordant links are concordant as links in
$\Z_p$-homology spheres via a concordance obtained by a similar
covering construction. Using the Blanchfield form of covering links,
they proved the long-standing conjecture that there are links which
are not concordant to boundary links in~\cite{CO1,CO2}.  Milnor's
$\bar\mu$-invariants~\cite{M2,M3} is also generalized for covering
links in~\cite{CO4}.

In this paper, we view covering links as links in rational homology
spheres, and utilize the signature invariant developed by the authors
in~\cite{CK3} to study covering links.  For homology boundary links,
we develop a new systematic method to compute signature jump functions
of covering links.  Recall that a link $L$ with $m$ components is
called a \emph{homology boundary link} if there exists an epimorphism
of the fundamental group of the complement of $L$ onto the free group
of rank~$m$.  An $m$-tuple $r=(r_1,\ldots,r_m)$ is called a
\emph{pattern} for $L$ if $r_i$ is the image of the $i$-th meridian
under the epimorphism~\cite{CL}.  A homology boundary link admits a
system of ``singular'' Seifert surfaces, and Seifert pairings and
Seifert matrices on singular Seifert surfaces are defined as
in~\cite{CO3}.  In Section~\ref{sec:covering-link}, we prove an
explicit formula to compute Seifert matrices and the signature jump
functions of covering links of a homology boundary link from its
Seifert matrix and pattern (see
Theorem~\ref{thm:covering-seifert-matrix-2comp}).

In order to prove the formula, we construct covering links and
their generalized Seifert surfaces in the sense of~\cite{CK3}
using cut-paste arguments, and compute Seifert matrices and
signature jump functions from the generalized Seifert surfaces.
The only algebraic tool needed is linear algebra of
matrices over the complex field. So our approach is geometric and
elementary in contrast to that in~\cite{CO3} where invariants for
covering links are obtained from the Blanchfield duality and their
invariance under concordance is shown using heavy machinery of
homological algebra.

There is another known way to compute Seifert matrices and signature jump
functions of covering links. For many concrete examples, the ambient space
of a covering link can be calculated as a surgery diagram using the method
of~\cite{AK}, and then we can construct a Seifert surface in the surgery
diagram and compute a Seifert matrix as illustrated in an example
of~\cite{CK3}.  Comparing with this, our approach is more practical in the
sense that we can compute invariants of a covering link of a given link
directly from the given link using a formula
without appealing to any diagram.

In Section~\ref{sec:app-to-link-conc}, we apply the above results on
covering links to study link concordance in spheres as influenced
by~\cite{CO1,CO2}. First we study links which are not concordant to
boundary links using the method of Section~\ref{sec:covering-link}.
Recall that a link is called a \emph{boundary link} if its components
bound disjoint Seifert surfaces.  Because boundary link concordance
classes can be algebraically classified~\cite{CS,Ko,Mi}, it had been
the center of interest whether all $(2q-1)$-links (with vanishing
$\bar\mu$-invariants if $q=1$) are concordant to boundary links.  The
first counterexamples was given by Cochran and Orr~\cite{CO1,CO2} as
mentioned before.  Gilmer and Livingston~\cite{GL}, and
Levine~\cite{L4} showed the same result using different techniques
like Casson-Gordon invariants and $\eta$-invariants, respectively. We
offer another method to detect links not concordant to boundary links
using signature jump functions of covering links. Roughly speaking, it
is shown that signature jump functions of covering links of links
which are concordant to boundary links must have period~$2\pi$.  Using
this, we show that if a homology boundary link has a specific form of
pattern and Seifert matrix, then it is not concordant to boundary
links (see Theorem~\ref{thm:link-not-concordant-to-boundary-link}).

Furthermore we show that there are 1-dimensional links with vanishing
$\bar\mu$-invariants which are positive mutants of ribbon links but
not concordant to boundary links (see
Theorem~\ref{thm:mutant-not-conc-to-boundary-link}).  For a link $L$
and a 3-ball $B$ in $S^3$ such that $L$ and $\partial B$ transversally
meet at exactly 4 points, the link obtained by pasting $(B,L\cap B)$
and $(S^3-\inte B, L-\inte B)$ along an orientation preserving
involution on $(\partial B,L\cap \partial B)$ whose fixed points are
disjoint from $L\cap \partial B$ is called a \emph{mutant} of~$L$. If
$L$ is oriented and the orientation of $L$ is preserved by the
mutation, then it is called a \emph{positive mutant}. Many link
invariants fail to distinguish links from their mutants.  It is known
that mutation preserves link invariants like Alexander, Jones,
Kauffman and HOMFLY polynomials, and positive mutation preserves link
signatures and $S$-equivalence classes of knots.  The problem to
distinguish links from positive mutants \emph{up to concordance} is
even subtler. The only known result is that Casson-Gordon invariants
are effective to distinguish some knots from their positive mutants up
to concordance, due to Kirk and Livingston~\cite{KL}. Almost nothing
has been known about the effect by (positive) mutation on link
concordance classes beyond knot concordance.  Our result says that
both of the set of slice links and the set of links concordant to
boundary links are not closed under positive mutations. We remark that
our result can also be viewed as a generalization of the result
of~\cite{JK} where it was shown that there is a boundary link with a
mutant which is not a (homology) boundary link. However since both the
link and its mutant in~\cite{JK} are ribbon links, it says nothing up
to concordance.

Following techniques of the classification of boundary link
concordance classes using Seifert matrices~\cite{Ko}, the appropriate
concordance classes of homology boundary links with a given pattern
are classified by Cochran and Orr~\cite{CO3}.  Thus it is now more
natural to ask whether all links (with vanishing $\bar\mu$-invariants
if $q=1$) are concordant to homology boundary links, instead of
boundary links.  As a partial answer, we show the following result by
investigating signatures of covering links of homology boundary links.

\begin{thm}\label{thm:link-not-concordant-to-hbl-pattern}
For any finite collection of patterns, there exist infinitely many
homology boundary links which are never concordant to any homology
boundary link admitting a pattern in the given collection.
\end{thm}

Note that every pattern is realized by a ribbon link~\cite{CL}. Combined
with our result, it can be seen that the variety of patterns which arise
in a concordance class of a homology boundary link depends heavily on the
choice of the concordance classes.

\section{Seifert matrices of covering links}\label{sec:covering-link}

In this section we derive formulae to compute Seifert matrices and
signature jump functions (of unions of parallel copies of components)
of covering links of homology boundary links.  It seems natural to
expect such formulae.  Firstly, Seifert matrices together with
patterns have enough information to classify the appropriate
concordance classes of homology boundary links~\cite{CO3}. Since the
signature jump function is invariant under link
concordance~\cite{CK3}, it is expected that signatures can be
calculated from Seifert matrices and patterns.  Secondly, since the
Blanchfield form is determined by a Seifert matrix~\cite{Ke,Hi,CO3}
and the Blanchfield form of a covering link of a link $L$ is the image
of the Blanchfield form of $L$ under a transfer
homomorphism~\cite{CO2,CO3}, it is expected that a Seifert matrix of a
covering link of $L$ can also be obtained from a Seifert matrix of
$L$. In this sense, our formula for Seifert matrices is analogous to
the transfer homomorphism for Blanchfield forms. We remark that no
explicit formula for the latter is known.

Throughout this paper, we consider ordered and oriented links only. We
use the following notations introduced in~\cite{CK3} for parallel
copies.  For a framed submanifold $M$ in an ambient space and an
$n$-tuple $\alpha=(s_1,\ldots,s_n)$ with $s_i=\pm 1$, let $i_\alpha M$
be the union of $n$ parallel copies of $M$, where the $i$-th copy is
oriented according to the sign of $s_i$, and let $n_\alpha$ be the sum
of~$s_i$. For a nonzero integer $r$, let $i_rM$ be the union of $|r|$
parallel copies of $M$ oriented according to the sign of~$r$.

We will consider only two component links to simplify notations,
though the arguments of this section can also be applied for links
with more than two components. Suppose that $L$ is a
$(2q-1)$-dimensional homology boundary link with components $J$ and
$K$ in a $\Z_p$-homology sphere $\Sigma$.  Let $\{E,F\}$ be a system
of singular Seifert surfaces properly embedded in the exterior $E_L$
such that $\partial E$, $\partial F$ are homologous to $J$, $K$ in a
tubular neighborhood of $L$, respectively.  An epimorphism from
$\pi_1(E_L)$ onto the free group on $x$ and $y$ is obtained by a
Thom-Pontryagin construction on $E\cup F$. By choosing meridians $\nu$
and $\mu$ based at a fixed basepoint outside $E\cup F$, a pattern
$r=(v,w)$ is determined which satisfies $v\equiv x$ and $w \equiv y$
modulo commutators.

Let $p$ be a prime and $d=p^a$ for some positive integer~$a$. Let
$\tilde\Sigma$ be the $d$-fold cyclic cover of $\Sigma$ branched
along~$J$, and $t$ be a generator of covering transformations.  Fixing
a basepoint of~$\tilde \Sigma$, the lift of $\mu$ based at the
basepoint is a meridian of a component $\tilde K$ of the pre-image
of~$K$. Then the union $\bigcup_{k=0}^{d-1} t^k \tilde K$ is a
covering link of~$L$. We consider the link $\tilde L = \bigcup i_{r_k}
t^k\tilde K$, where the parallel copies are taken with respect to the
framing induced by $E$ and~$F$.  We will compute the signature jump
function $\delta_{\tilde L}(\theta)$ defined in~\cite{CK3} from the
given data $r$ and~$\{E, F\}$.

We will construct $\tilde \Sigma$ using well-known cut and paste
arguments as in~\cite{AK,CK2,Hi,K}, and construct a Seifert surface of
$\tilde L$ by taking parallel copies of lifts of~$F$.  Denote the
pre-image of $J$ in $\tilde\Sigma$ by~$\tilde J$. Choose a smaller
tubular neighborhood $V$ of $L$ in~$\Sigma-E_L$.  We can cancel out
boundary components of $E$ with opposite orientations by attaching to
$E$ annuli properly embedded in~$\Sigma-\inte (E_L\cup V)$, and then
we obtain a proper submanifold $N$ in $\Sigma-\inte V$ such that
$\partial N$ is a single parallel of $J$ on $\partial V$.  Removing
from $\Sigma$ the interior of the component of $V$ containing $J$, we
obtain an exterior $E_J$.  Choose a bicollar $N\times[-1,1]$ in $E_J$
so that $N\times 1$ is a translation of $N$ along the positive normal
direction. For $k=0,\ldots,d-1$, let $t^k \tilde X$ be a copy of
$X=E_J-N\times(-1,1)$ and $g^k_{\pm}\colon N \to t^k\tilde X$ be a
copy of the inclusions $g_\pm\colon N \to N\times\{\pm 1 \}\subset X$.
Then the exterior $E_{\tilde J}$ of $\tilde J$ in $\tilde\Sigma$ is
homeomorphic to the quotient space
$$
\Big(\bigcup_{k=0}^{d-1} t^k\tilde X \Big)\Big/ \sim
$$
where $g^k_+(z)$ and $g^{k+1}_-(z)$ (indices are modulo $d$) are
identified for $z\in N$. $\tilde\Sigma$ is obtained by gluing
$E_{\tilde J}$ and $S^q\times D^2$ along boundaries.  We remark that
we can construct $E_{\tilde J}$ using $E$ instead of $N$. The reason
why we use $N$ is that there is a duality isomorphism between
$H_q(X;\Q)\cong H_q(E_J-N;\Q)$ and $H_q(N;\Q)$ for any $q$. This
isomorphism will be needed later and is not established for $q=1$ if
we use $E$ instead of $N$.

Let $t^k \tilde N$ = $g^k_+(N) \subset E_{\tilde J}$, and
denote the lift of $F$ in $t^k\tilde X$ by $t^k \tilde F$.  Then
$$
\partial(t^k\tilde F) = \Big(\bigcup_{l=0}^{d-1}
i_{\alpha_{kl}} t^l\tilde K \Big) \cup i_\alpha \tilde J
$$
for some tuples $\alpha$ and $\alpha_{kl}$.  Obviously $n_\alpha=0$.
$n_{\alpha_{kl}}$ is determined by the pattern as follows.  Since $w
\equiv y$ modulo commutators, we can write $w = \prod_i x^{a_i}y^{b_i}
x^{-a_i}$ where $b_i=\pm 1$ and $\sum_i b_i = 1$. Let $c_n(r)$ be the
sum of $b_i$ over all $i$ such that $a_i=n$.  If we travel along the
lift of $\mu$ which is a meridian of $\tilde K$ in $\tilde\Sigma$, a
$\pm$-intersection with $t^{a_i} \tilde F$ occurs for each
$x^{a_i}y^{\pm 1}x^{-a_i}$ factor in $w$.  From this observation,
$n_{\alpha_{kl}}$ is the sum of $c_n(r)$ over all $n$ satisfying
$n\equiv k-l$ mod $d$.  We remark that for any pattern $r$, all but
finitely many $c_n(r)$ vanish, and $\sum_n c_n(r) = 1$.

The following lemma implies that for any $r_0,\ldots,r_{d-1}$, the
system of $d$ equations
$$
\sum_{k=0}^{d-1}n_{\alpha_{kl}} x_k= r_l \quad(l=0,\ldots,d-1)
$$
has a unique solution $(x_k)$ over $\Q$.

\begin{lem}\label{lem:complexity-matrix-det}
If $n$ is a prime power and $c_1,\ldots,c_n$ are integers such that
$c_1+\cdots+c_n=1$,
$$
\begin{bmatrix}
c_1 & c_2 & \cdots & c_{n} \\
c_{n} & c_1 & \cdots & c_{n-1} \\
\vdots & \vdots & \ddots & \vdots \\
c_2 & c_3 & \cdots & c_1
\end{bmatrix}
$$
is a nonsingular matrix.
\end{lem}
\begin{proof}
Let $n=p^a$ and $p$ be a prime. First we observe that the matrix has a
symmetry in the sense that it is invariant under the $\Z_n$-action which
cyclically shifts rows and columns.

We expand the determinant as a sum over all permutations of
$\{1,\ldots,n\}$, and investigate when a particular monomial, to say,
$m=c_1^a c_{i_1}^{a_1}\cdots c_{i_k}^{a_k}$ $(1<i_1<\cdots<i_k)$
appears as a summand.  Let $X$ be the set of all subsets of
$\{1,\ldots,n\}$.  The action on $\{1,\ldots,n\}$ by $\Z_n$ induces on
$X$ in an obvious way.  For any element $x$ of $X$ with
cardinality~$a$, let $P_x$ be the set of permutations $\pi$ such that
$x$ is the fixed point set of $\pi$ and the product of the
$(i,\pi(i))$-th entries is equal to~$m$.  If two elements $x$ and $y$
in $X$ are in the same orbit, then the action induces a bijection
between $P_x$ and $P_y$ which preserves the signs of permutations, by
the symmetry.  Hence the coefficient of $m$ in the determinant is an
integral linear combination of the cardinalities of orbits. If
$0<a<n$, the cardinality of an orbit is a multiple of~$p$, and so is
the coefficient of $m$.  By the symmetry again, the same argument
works when we replace $c_1$ by any~$c_i$, and this shows that
coefficients of all monomials except $c_i^n$ are multiples of~$p$.
Therefore the determinant is congruent to $c_1^n+\cdots+ c_n^n \equiv
c_1+\cdots+c_n \equiv 1$ modulo~$p$.
\end{proof}

Let $s$ be a common multiple of denominators of $x_k$ and let
$$
M=\bigcup_{k=0}^{p-1} i_{sx_k} t^k\tilde F.
$$
Then we have
$$
\partial M =
\Big(\bigcup_{k,l} i_{sx_k} i_{\alpha_{kl}} t^l\tilde K \Big)
\cup \Big(\bigcup_k i_{sx_k} i_{\alpha_k} \tilde J \Big) =
\Big(\bigcup_l i_{\beta_l} t^l\tilde K\Big)
\cup i_\beta \tilde J
$$
where $\beta_l$, $\beta$ are tuples such that $n_{\beta_l} = \sum_k
sx_k n_{\alpha_{kl}} = sr_l$, $n_\beta = 0$.  By attaching annuli to
$M$ in a tubular neighborhood of $(\bigcup t^k \tilde K)\cup \tilde J$
to cancel out unnecessary boundary components, we obtain a submanifold
$M'$ with boundary $i_s\tilde L$ and we can compute $\delta_{\tilde
L}(\theta)$ from a Seifert matrix of $M'$.  For $q>1$, $\delta_{\tilde
L}(\theta)$ can be computed from a Seifert matrix $P$ of $M$ since
$H_q(M) \cong H_q(M')$.

For $q=1$, we need additional arguments.  Let $S\colon H_1(M')\times
H_1(M')\to \Q$ be the Seifert pairing of $M'$.  For a manifold $V$,
denote the cokernel of $H_i(\partial V)\to H_i(V)$ by $\bar H_i(V)$.
Then we have $H_1(M')\cong \bar H_1(M)\oplus\Z^{2n}\oplus\Z^m$, where
the $\Z^{2n}$ factor is generated by cores of attached annuli and
their dual loops, and the $\Z^m$ factor is generated by boundary
parallel loops. We will show that $S$ induces a well-defined ``Seifert
pairing'' on $\bar H_1(M)$ and the $\Z^{2n}\oplus \Z^m$ factor has no
contribution to the signature \emph{jump} function of $S$.  Thus
$\delta_{\tilde L}(\theta)$ can be computed from a Seifert matrix $P$
defined on $\bar H_1(M)$. Hence we can unify notations for any $q$ by
letting $P$ be a Seifert matrix on $\bar H_q(M)$ and we have
$\delta_{\tilde L}(\theta) = \delta^q_P(\theta/s)$.

Our assertion for $q=1$ is shown as follows.  We choose generators
$\{c_i, d_i\}$ and $\{e_i\}$ of $\Z^{2n}$ and $\Z^m$ factor,
respectively, where $c_i$ is the core of an attached annulus, $d_i$ is
a curve on $M'$ whose intersection number with $c_j$ is $\delta_{ij}$
(Kronecker's delta symbol), and $e_i$ is a boundary component of
$M$. We will show the linking number of a loop on $M$ and a boundary
component $c$ of $M$ is zero. We may assume $c=t^j\tilde K$ or $\tilde
J$ since $c$ is homologous to one of them in $\tilde\Sigma-\inte
M$. For any $j=0,\ldots,d-1$, the equation $\sum_k
n_{\alpha_{kl}}x_k=\delta_{jl}$ has a solution by
Lemma~\ref{lem:complexity-matrix-det}, and so we can construct a
surface in $\tilde\Sigma$ whose boundary is homologous to $i_a
t^j\tilde K$ for some $a>0$ by taking parallel copies of $t^k \tilde
F$ as before. By attaching annuli, we obtain a surface that is
disjoint to $M$ and bounded by $i_a t^j \tilde K$.  Therefore the
linking number of $t^j \tilde K$ and any loop on $M$ is zero.
Similarly the linking number of $\tilde J$ and any loop on $M$ is zero
since we can construct a surface which is bounded by $\tilde J$ and
disjoint to $M$ by attaching annuli to $\tilde N$. Since $c_i$ and
$e_i$ are homologous to boundary components of $M$ and $M$ induces
0-linking framings of boundary components, the Seifert pairing $S$
vanishes on the pairs $(c_i,x)$, $(x,c_i)$, $(e_i,x)$, $(x,e_i)$,
$(c_i,e_j)$, $(e_i,c_j)$, $(c_i, c_j)$, and $(e_i,e_j)$ for any $x$ in
$H_1(M)$.  By the choice of $c_i$ and $d_j$,
$S(c_i,d_i)-S(d_i,c_i)=\delta_{ij}$. From the observations, the usual
Seifert pairing determines a well-defined ``Seifert pairing'' on $\bar
H_1(M)$, and furthermore the Seifert matrix $Q$ over the chosen basis
of $H_1(M')$ is given by
$$
{\arraycolsep=.5em \def\arraystretch{1.5}
Q=
\left[\begin{array}{c|c|c|c}
P & 0 & *     & 0 \\
\hline
0 & 0 & R^T+I & 0 \\
\hline
* & \setbox0=\hbox{$\hphantom{R^T+I}$} \hbox to\wd0{\hss $R$\hss} & *     & * \\
\hline
0 & 0 & *     & 0 \\
\end{array}\right]
}
$$
where $P$ is a Seifert matrix defined on $\bar H_1(M)$, and $R$
represents the Seifert pairing between bases $\{d_i\}$ and $\{c_i\}$.
In order to compute $\sigma^+_{Q}(\phi)$, we consider a complex
hermitian matrix
$$
{\arraycolsep=.5em \def\arraystretch{2}
\displaystyle \frac{wQ-Q^T}{w-1} =
\left[\begin{array}{c|c|c|c}
\displaystyle \frac{wP-P^T}{w-1} & 0 & * & 0 \\[1ex]
\hline
0 & 0 & \displaystyle R^T+\frac{w}{w-1}I & 0 \\[1ex]
\hline
* & \displaystyle R^{\strut}-\frac{1}{w-1}I & * & * \\[1ex]
\hline
0 & 0 & * & 0
\end{array}\right]
}
$$
for an uni-modular complex number $w$.  The submatrix $R-(w-1)^{-1}I$
can be viewed as a matrix over the ring of polynomials in
$z=(w-1)^{-1}$. Since the determinants of the upper-left square
submatrices of $R-(w-1)^{-1}I$ are nonzero polynomials in $z$, the
pivots used in the Gauss-Jordan elimination process are nonzero
rational functions in $z$ and hence $R-(w-1)^{-1}I$ can be transformed
into a nonsingular diagonal matrix by row operations on
$(wQ-Q^T)/(w-1)$ if $w$ is not a zero of the denominators and the
numerators of the pivots. Since $(wQ-Q^T)/(w-1)$ is hermitian, the
submatrix $R^T+w(w-1)^{-1}I$ is transformed into a nonsingular
diagonal matrix by corresponding column operations.  Note that this
can be performed all but finitely many $w$, and it does not alter
vanishing blocks of $(wQ-Q^T)/(w-1)$. By further row and column
operations on $(wQ-Q^T)/(w-1)$, all outer blocks except the top-left
block are cleared and we eventually obtain the block sum of
$(wP-P^T)/(w-1)$, a nonsingular null-cobordant matrix and a zero
matrix.  Therefore $\sigma^+_Q(\phi)=\sigma^+_P(\phi)$ on a dense
subset of $\R$.  This shows our assertion for $q=1$.

Now we need to compute the Seifert matrix $P$ defined on $\bar
H_q(M;\Q)$.  (Note that if $q>1$, $\bar H_q(-)$ is identified with
$H_q(-)$ and $P$ is the usual Seifert matrix of $M$.) $P$ is obtained
by duplicating rows and columns of a Seifert matrix defined on $\bar
H_q(\bigcup t^k\tilde F;\Q)\cong \bigoplus\bar H_q(F;\Q)$, which we
will compute.

Let $x$ and $y$ be elements of $\bar H_q(F;\Q)$ and let $a$ and $b$ be
$q$-cycles on $F$ which represent the image of $x$ and $y$ under a
fixed splitting map $\phi\colon \bar H_q(F;\Q)\to H_q(F;\Q)$,
respectively.  We will compute the linking number of lifts $t^k\tilde
a^+$ and $\tilde b$ in $\tilde\Sigma$, where $\tilde z$ denotes the
unique lift of $z$ in $\tilde X$ for a chain $z$ in $X$. A Seifert
pairing on $\bar H_q(N\cup F;\Q)\cong H_q(N;\Q)\oplus \bar H_q(F;\Q)$
is induced by $\phi$ and the usual Seifert pairing on $H_q(N\cup F)$.
Fix basis of $H_q(N;\Q)$ and $\bar H_q(F;\Q)$, and let $\spmatrix{A&B
\\ \epsilon B^T & C}$ be a Seifert matrix defined on $H_q(N;\Q)\oplus
\bar H_q(F;\Q)$ with respect to the basis as in~\cite{Ko,CO3}.  By
duality, we have $H_q(N;\Q) \cong H_q(X;\Q)$.  $(g_+)_*, (g_-)_*\colon
H_q(N;\Q)\to H_q(X;\Q)$ and the composition of $\phi$ and $g_* \colon
H_q(F;\Q)\to H_q(X;\Q)$ are represented by $A$, $\epsilon A^T$ and
$B$, respectively.

We will find $q$-cycles $z_0,\ldots,z_{d-1}$ on $N$ and $(q+1)$-chains
$u_1,\ldots,u_{d-1}$ in $X$ such that
\begin{align*}
g(b) + g_+(z_0) - g_-(z_1) &= \partial u_0 \\
       g_+(z_1) - g_-(z_2) &= \partial u_1 \\
                           &\vdots         \\
       g_+(z_{d-2}) - g_-(z_{d-1}) &= \partial u_{d-2} \\
       g_+(z_{d-1}) - g_-(z_0) &= \partial u_{d-1}.
\end{align*}
Once finding $z_i$ and $u_i$, we obtain a chain $\tilde u = \bigcup
t^k \tilde u_k$ in $\tilde\Sigma$ such that $\partial\tilde u = \tilde
b$. Then we can compute the linking number of lifts of $a$ and $b$ as
follows.
\begin{equation*}\begin{split}
\lk_{\tilde\Sigma}(t^k\tilde a^+,\tilde b) &=
t^k\tilde a^+\cdot\tilde u = a^+ \cdot u_k\\
&= \lk_{\Sigma}(a^+,\partial u_k) \\
&= \begin{cases}
\lk_{\Sigma}(a^+,b)+\lk_{\Sigma}(a^+,z_0)-\lk_{\Sigma}(a^+,z_1), & k=0, \\
\lk_{\Sigma}(a^+,z_k)-\lk_{\Sigma}(a^+,z_{k+1}), & 1\le k \le d-1.
\end{cases}
\end{split}\end{equation*}

Viewing $x$, $y$, $z_k$ as column vectors representing elements of
appropriate $\Q$-homology groups, the above system of equations
becomes
$$
\begin{bmatrix}
A & -\epsilon A^T & & \\
& A & -\epsilon A^T & \\
& & \ddots & \ddots \\
& & & A & -\epsilon A^T \\
-\epsilon A^T & & & & A
\end{bmatrix}
\begin{bmatrix}
z_0 \\
z_1 \\
\vdots \\
z_{d-2} \\
z_{d-1}
\end{bmatrix} =
\begin{bmatrix}
-B y \\
0 \\
\vdots \\
0 \\
0
\end{bmatrix}
$$
in $H_q(X;\Q)$. Since $N$ has one boundary component, $A-\epsilon
A^T$ is nonsingular. By multiplying $(A-\epsilon A^T)^{-1}$ on the
left of each row, it becomes
$$
\begin{bmatrix}
\Gamma & I-\Gamma & & \\
& \Gamma & I-\Gamma & \\
& & \ddots & \ddots \\
& & & \Gamma & I-\Gamma \\
I-\Gamma & & & & \Gamma
\end{bmatrix}
\begin{bmatrix}
z_0 \\
z_1 \\
\vdots \\
z_{d-2} \\
z_{d-1}
\end{bmatrix} =
\begin{bmatrix}
-(A-\epsilon A^T)^{-1} By \\
0 \\
\vdots \\
0 \\
0
\end{bmatrix}
$$
where $\Gamma=(A-\epsilon A^T)^{-1}A$ and $I$ is the identity matrix.
Since $\Gamma^d-(\Gamma-I)^d$ is a presentation matrix of $H_q(\tilde
\Sigma)$ (see \cite{Se}) and $\tilde \Sigma$ is a rational homology
sphere, $\Gamma^d-(\Gamma-I)^d$ is invertible and a unique solution
$(z_k)$ exists.  It is easy to check that
$$
z_k = \begin{cases}
\displaystyle
-\frac{\Gamma^{d-1}}{\Gamma^d-(\Gamma-I)^p}
       (A-\epsilon A^T)^{-1} By, & k=0, \\[2ex]
\displaystyle
-\frac{\Gamma^{k-1}(\Gamma-I)^{d-k}}{\Gamma^d-(\Gamma-I)^d}
       (A-\epsilon A^T)^{-1} By, & 1 \le k \le d-1. \\
\end{cases}.
$$
Note that in the above fractional notations of matrices,
denominators and numerators commute and so we have no ambiguity.
By the above calculation of the linking number, we have
$$
\lk_{\tilde\Sigma}(t^k\tilde a^+, \tilde b) =
\begin{cases}
\displaystyle x^T \Big(C-\epsilon B^T
\frac{\Gamma^{d-1}-(\Gamma-I)^{d-1}}{\Gamma^d-(\Gamma-I)^d}
(A-\epsilon A^T)^{-1} B\Big)y,
& k=0, \\[2ex]
\displaystyle x^T \Big(\epsilon B^T
\frac{\Gamma^{k-1}(\Gamma-I)^{d-k-1}}{\Gamma^d-(\Gamma-I)^d}
(A-\epsilon A^T)^{-1} B\Big)y,
& 1 \le k \le d-1.
\end{cases}
$$

From the above discussion, we obtain the following result.

\begin{thm}\label{thm:covering-seifert-matrix-2comp}
Let $\spmatrix{A & B\\ B^T & C}$ be a Seifert matrix of $L$ defined on
$E\cup F$ in the sense of~\cite{CO3}. Then the block matrix
$(A_{kl})_{0\le k,l<d}$ given by
$$
A_{kl}=
\begin{cases}
\displaystyle C-\epsilon B^T
\frac{\Gamma^{p-1}-(\Gamma-I)^{p-1}}{\Gamma^p-(\Gamma-I)^p}
(A-\epsilon A^T)^{-1} B,
& k=l, \\[2ex]
\displaystyle \epsilon B^T
\frac{\Gamma^{k-l-1}(\Gamma-I)^{p-k+l-1}}{\Gamma^p-(\Gamma-I)^p}
(A-\epsilon A^T)^{-1} B,
& k > l,\\[2ex]
\displaystyle \epsilon B^T
\frac{\Gamma^{p-k+l-1}(\Gamma-I)^{k-l-1}}{\Gamma^p-(\Gamma-I)^p}
(A-\epsilon A^T)^{-1} B,
& k < l,
\end{cases}
$$
is (cobordant to if $q=1$) a Seifert matrix defined on $\bar H_q(\bigcup
t^k\tilde F;\Q)$, and the block matrix $(P_{kl})_{0\le k,l<d}$ given by
$$
P_{kl} =
\begin{cases}
i^q_{sx_k} A_{kl}, & k=l \\
\text{\rm $sx_k\times sx_l$ array of $A_{kl}$}, & k\ne l
\end{cases}
$$
is (cobordant to if $q=1$) a Seifert matrix defined on $\bar H_q(M;\Q)$.

In particular, $\delta_{\tilde L}(\theta)=\delta^q_{(P_{kl})}(\theta/s)$.
\end{thm}

\begin{proof}
For $q>1$, we have already proved the theorem. For $q=1$, we have
proved that the conclusion holds if $\spmatrix{A & B\\ B^T & C}$ is a
Seifert matrix defined on $H_1(N)\oplus \bar H_1(F)$.  By observing
that the first formula induces a well-defined homomorphism
$G(2,\epsilon) \to G(d,\epsilon)$ which sends $\spmatrix{A & B\\ B^T &
C}$ to $(A_{kl})$ on the groups of cobordism classes of Seifert
matrices in the sense of~\cite{Ko}, it suffices to show that a Seifert
matrix defined on $\bar H_1(E)\oplus \bar H_1(F)$ in~\cite{CO3} and a
Seifert matrix defined on $H_1(N)\oplus \bar H_1(F)$ in the previous
discussion are cobordant in the sense of~\cite{Ko} since both Seifert
matrices represent elements of $G(2,\epsilon)$.

This assertion is proved by a similar reduction argument used earlier
for Seifert matrices on $H_1(M')$ and $\bar H_1(M)$.  $H_1(N)\cong
\bar H_1(E)\oplus \Z^{2n}$ where the $\Z^{2n}$ factor is generated by
cores of annular components of $N-\inte E$ and its dual generators,
and the linking number of each core and any loop on $E\cup F$ is zero.
Therefore
$$
{\arrayrulewidth=.2pt \doublerulesep=\arrayrulewidth
\left[\begin{array}
{c|ccccc||c}
A & * & 0 & \cdots & * & 0 & B \\[.2ex]
\hline
* & * & * & \cdots & * & * & * \\
0 & * & 0 & \cdots & * & 0 & 0 \\
\vdots & \vdots & \vdots & \ddots & \vdots & \vdots & \vdots \\
* & * & * & \cdots & * & * & * \\
0 & * & 0 & \cdots & * & 0 & 0 \\
\hline\hline
B^T\vphantom{B^{T^T}} & * & 0 & \cdots & * & 0 & C
\end{array}\right]
}
$$
is a Seifert matrix defined on $H_1(N)\oplus \bar H_1(F)$, where
$\spmatrix{A & B\\ B^T & C}$ is a Seifert matrix defined on $\bar
H_1(E)\oplus \bar H_1(F)$. It is easy to check that the block sum of
this Seifert matrix and $-\spmatrix{A & B\\ B^T & C}$ is
null-cobordant.
\end{proof}

\section{Application to link concordance}
\label{sec:app-to-link-conc}

\subsection*{Concordance of boundary links}
In this subsection we study examples of homology boundary links in
$S^{2q-1}$ which are not concordant to boundary links, whose existence
was shown first in~\cite{CO1,CO2} and subsequently in \cite{GL,L4}.  A
key observation in~\cite{CO2} is that a covering link of a boundary
link $L$ is again a boundary link. Since a boundary link is a
primitive link in the sense of~\cite{CK3}, the signature jump function
of (any union of parallels of components of) a covering link of $L$
must be of period $2\pi$ by Theorem~1.2 of~\cite{CK3}.  Since
corresponding covering links of concordant links are concordant, the
same conclusion holds under an weaker assumption that $L$ is
concordant to a boundary link by the fact that signatures are
invariants under link concordance~\cite{CK3}.  We state this as a
theorem.

\begin{thm}\label{thm:covering-link-of-boundary-link}
If a link $L$ is concordant to a boundary link, the signature jump
function of any union of parallels (with respect to the 0-linking
framing if $L$ is 1-dimensional) of components of a covering link of
$L$ has the period $2\pi$.
\end{thm}

Using Theorem~\ref{thm:covering-link-of-boundary-link}, we prove

\begin{thm}\label{thm:link-not-concordant-to-boundary-link}
Suppose $L$ is a 2-component homology boundary link in $S^{2q+1}$ with
a pattern $r$ and a Seifert matrix $\spmatrix{A & B \\ \epsilon B^T &
C}$ in the sense of~\cite{CO3} such that $A = C = \spmatrix{V & V \\
\epsilon V^T & \epsilon V^T}$, $B = \spmatrix{ V & V \\ \epsilon V^T &
V}$ for a Seifert matrix $V$ of a knot with nontrivial signature jump
function, and for some $n_0$, $c_n(r)=0$ if and only if $n \ne
n_0,n_0+1$.  Then $L$ is not concordant to any boundary links.
\end{thm}

\begin{proof}
In this proof, we denote components of $L$ by $J$, $K$, and use the
notations of Section~\ref{sec:covering-link}.  We consider the
covering link of $L$ obtained by taking the $p$-fold cyclic cover of
$S^{2q+1}$ branched along the first component $J$ for an odd
prime~$p$.  Let denote the first component of the pre-image of $K$
by~$\tilde K_L$.

We have $\Gamma =
(A-\epsilon A^T)^{-1}A = \spmatrix{ G & G \\ 1-G & 1-G}$ where $G =
(V-\epsilon V^T)^{-1}V$. By a straightforward calculation using
Theorem~\ref{thm:covering-seifert-matrix-2comp} and the fact
$\Gamma^2=\Gamma$,
$$
(A_{kl})=\left[ \begin{array}{cc|cc|ccc|cc}
V & & & & & & & & V \\
 & \epsilon V^T & \epsilon V^T & & & & & & \\
\cline{1-5} \cline{7-9}
 & V & V & & & & & & \\
 & & & \epsilon V^T & \epsilon V^T & & & & \\
\cline{1-5} \cline{7-9}
 & & & V & V & & & & \\
\multicolumn{5}{c}{} & \ddots & \multicolumn{2}{c}{} \\
 & & & & & & \epsilon V^T & \epsilon V^T & \\
\cline{1-5} \cline{7-9}
 & & & & & & V & V & \\
\epsilon V^T  & & & & & & & & \epsilon V^T
\end{array} \right]
$$
is a Seifert matrix defined on $\bar H_q(\bigcup_k t^k \tilde F;\Q)$.

Denote $c_n(r)$ by $c_n$ for simplicity.  By conjugating the pattern
by $x^{-n_0}$, we may assume $n_0=0$ and $c_0=m$, $c_1=1-m$ for some
$m\ne 0,1$.  Moreover by reversing orientations if necessary, we may
assume that $m>1$.  Since $x_0 = m^{p-1}/(m^p-(m-1)^p)$, $x_k =
m^{k-1}(m-1)^{p-k}/(m^p-(m-1)^p)$ $(k=1,\ldots,p-1)$ form a solution
of the linear system
$$
\begin{bmatrix}
c_0 & c_1 & \cdots & c_{p-1} \\
c_{p-1} & c_0 & \cdots & c_{p-2} \\
\vdots & \vdots & \ddots & \vdots \\
c_1 & c_2 & \cdots & c_0
\end{bmatrix}
\begin{bmatrix} x_0 \\ x_1 \\ \vdots \\ x_{p-1} \end{bmatrix}
=
\begin{bmatrix} 1\\ 0 \\ \vdots \\ 0 \end{bmatrix},
$$
We can compute the matrix $(P_{kl})$ in
Theorem~\ref{thm:covering-seifert-matrix-2comp} by putting
$s=m^p-(m-1)^p$.  $(P_{kl})$ is transformed to the block sum of $
i^q_{x_1-x_0}V, \ldots, i^q_{x_{p-1}-x_{p-2}}V,i^q_{x_0-x_{p-1}}V$ by
permuting rows and columns.  Therefore, by the reparametrization
formula in~\cite{CK3} and
Theorem~\ref{thm:covering-seifert-matrix-2comp}, we have
$$
\delta_{\tilde K_L}(\theta) = \sum_{k=0}^{p-1} \delta^q_V(y_k\theta)
$$
where $y_0 = -\frac{1-a^{p-1}}{m(1-a^p)}$, $y_k =
\frac{a^{k-1}}{m^2(1-a^p)}$ for $k=1,\ldots,p-1$ and $a=(m-1)/m$.

Since $\delta^q_V$ is nontrivial, there exists $\theta_0>0$ such that
$\delta^q_V(\theta_0)\ne 0$ and $\delta^q_V(\theta)= 0$ for all
$|\theta| < \theta_0$.  Let $\theta_1=\theta_0/y_0$. Then
$\delta_{\tilde K_L} (\theta_1) = \delta^q_V(\theta_0) \ne 0$ for
sufficiently large $p$, since for $k>0$ $|y_k/y_0|$ uniformly
converges to $1/m$ as $p\to\infty$.  We can choose large $N$ such that
$|y_k(\theta_0/y_0-2\pi)| < \theta_0$ for all $p,k>N$, since
$y_k(\theta_0/y_0-2\pi)$ uniformly converges to 0 as $k\to\infty$.
Therefore if $p>N$, $\delta_{\tilde K_L}(\theta_1-2\pi) =
\delta^q_V(\theta_0-2\pi y_0) + \delta^q_V(y_1(\theta_0/y_0-2\pi)) +
\cdots + \delta^q_V(y_N(\theta_0/y_0-2\pi))$.  Since
$\{y_k(\theta_0/y_0-2\pi)\}_{p=1}^\infty$ is a monotone convergent
sequence for each $k=0,\ldots,N$ and the set of points at which
$\delta^q_V$ is nonzero is discrete (see Lemma~2.1 of~\cite{CK3}),
$\delta^q_V(y_k(\theta_0/y_0-2\pi))=0$ for any large $p$. This shows
that $\delta_{\tilde K_L}(\theta)$ is not of period $2\pi$ for any
large $p$.  By Theorem~\ref{thm:covering-link-of-boundary-link}, $L$
is not concordant to boundary links.
\end{proof}

We remark that any knot that is not (algebraically if $q=1$) torsion
in the knot concordance group has a Seifert matrix $V$ satisfying the
hypothesis of the above theorem~\cite{L1,L2}.

\begin{figure}[hbt]
\begin{center}
\includegraphics{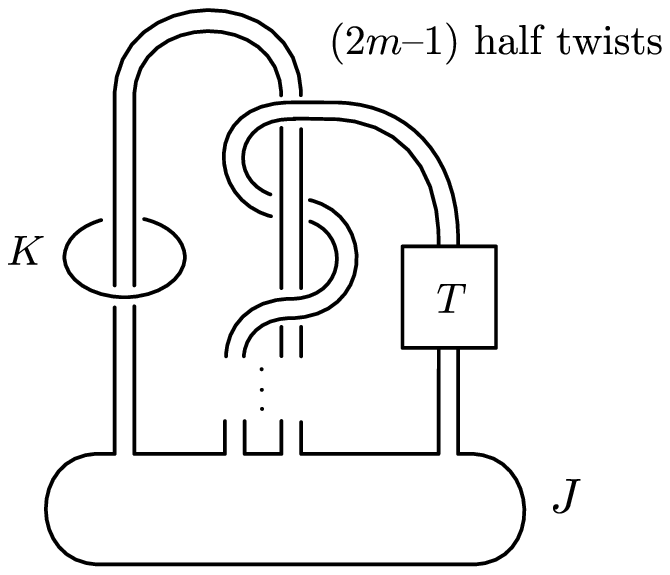}
\caption{}\label{fig:colink}
\end{center}
\end{figure}

In~\cite{CO3}, it is shown that an arbitrary pair of a pattern and a
Seifert matrix is always realized by a geometric construction of a
homology boundary link. Hence we can obtain a large collection of
links that satisfy the conditions of
Theorem~\ref{thm:link-not-concordant-to-boundary-link} and therefore
are not concordant to boundary links.  We remark that the main
examples of links in~\cite{CO2}, denoted by $L(T,m)$, also satisfy the
conditions of Theorem~\ref{thm:link-not-concordant-to-boundary-link}
for $m\ne 0,1$. In fact, the conditions of
Theorem~\ref{thm:link-not-concordant-to-boundary-link} can be viewed
as an algebraic description of $L(T,m)$. In the case of $q=1$,
$L(T,m)$ in the three space is illustrated in Figure~\ref{fig:colink}.
The first component $J$ bounds the obvious Seifert surface with one
0-handle and two 1-handles where a knot $T$ is tied along one of the
1-handles.  The Seifert matrix with respect to the generators
represented by the 1-handles is given by $\spmatrix{0 & m \\ m-1 &
0}$.

\subsection*{Mutation and link concordance}
\begin{figure}[hbt]
\begin{center}
\includegraphics{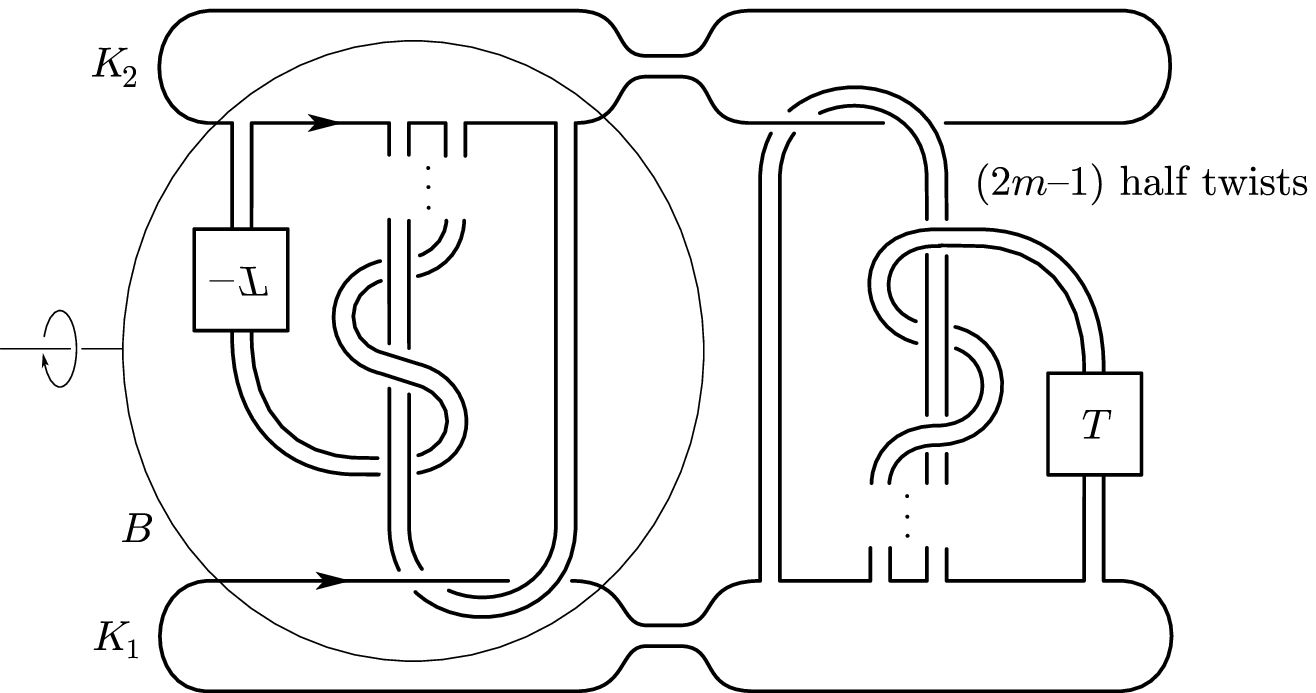}
\caption{}\label{fig:mutationex}
\end{center}
\end{figure}

In this subsection we illustrate an example of 1-dimensional links
which are positive mutants of ribbon links but not concordant to
boundary links.  Consider the 1-dimensional link $L$ with two
components shown in Figure~\ref{fig:mutationex}. The first component
$K_1$ of $L$ has the same knot type as that of~$L(T,m)$.  The other
component $K_2$ has the knot type of the mirror image of $K_1$.

\begin{thm}\label{thm:mutant-not-conc-to-boundary-link}
$L$ has the following properties:
\begin{enumerate}
\item
$L$ is a positive mutant of a ribbon link.
\item
$L$ is a fused boundary link. In particular, $L$ has vanishing
$\bar\mu$-invariants.
\item
If $T$ is not torsion in the algebraic knot concordance group and $m\ne
0,1$, $L$ is not concordant to boundary links.
\end{enumerate}
\end{thm}

\begin{proof}
By the positive mutation on the 3-ball $B$ shown in
Figure~\ref{fig:mutationex}, we obtain a mutant $L^*$ of~$L$. $L^*$ is
a connected sum of $L(T,m)$ and its mirror image, and in particular,
$L^*$ is a ribbon link.

$L$ is obtained by attaching two bands joining disjoint components of
the boundary link that is the split union of two parallel copies of
$T$ and $-T$ (the mirror image of $T$). Therefore $L$ is a fused
boundary link and has vanishing $\bar\mu$-invariants.

To show the last conclusion, we consider the first component $\tilde
K_L$ of the $p$-fold covering link of $L$ as before. $L$ is a
connected sum of $L(T,m)$ and $L'$, where $L'$ is the link obtained by
exchanging the order of the components of~$L(-T,-m)$. Hence $\tilde
K_L$ is the connected sum of $\tilde K_{L(T,m)}$ and $\tilde
K_{L'}$. By the additivity of signature jump function~\cite{CK3},
$\delta_{\tilde K_L}(\theta)=\delta_{\tilde
K_{L(T,m)}}(\theta)+\delta_{\tilde K_{L'}}(\theta)$. Since the first
component of $L'$ (which is the second component of $L(-T,-m)$) is
unknotted, the ambient space of $\tilde K_{L'}$ is the 3-sphere and
$\delta_{\tilde K_{L'}}(\theta)$ is of period~$2\pi$. (In fact,
$\tilde K_{L'}$ has the knot type of~$-(T \# T)$.) Therefore the
period of $\delta_{\tilde K_L}(\theta)$ is equal to that of
$\delta_{\tilde K_{L(T,m)}}(\theta)$, and is not equal to $2\pi$ for
any sufficiently large prime $p$ by the proof of
Theorem~\ref{thm:link-not-concordant-to-boundary-link}. This proves
that $L$ is not concordant to boundary links.
\end{proof}

\subsection*{Concordance of homology boundary links with given patterns}
In this subsection we generalize the previous arguments to show
Theorem~\ref{thm:link-not-concordant-to-hbl-pattern}.

\begin{proof}
[Proof of Theorem~\ref{thm:link-not-concordant-to-hbl-pattern}] Fix a
pattern $r$.  Suppose that $L=K_1\cup \cdots\cup K_n$ is a homology
boundary link admitting $r$ as a pattern. Let $F_1\cup\cdots\cup F_n$
be a singular Seifert surface.  Consider the covering link
$\bigcup_{i>1,k} t^k \tilde K_i$ of $L$ obtained by taking the
$p$-fold cover branched along $K_1$ as before.  Attaching annuli to a
lift of $F_2$ as done in Section~\ref{sec:covering-link}, we obtain a
submanifold in the ambient space of the covering link which is bounded
by $\tilde L = \bigcup_{i>1,k} i_{c_{ik}} t^k \tilde K_i$ for some
integers~$c_{ik}$. In particular, $\tilde L$ is a primitive link and
so $\delta_{\tilde L}(\theta)$ has the period $2\pi$.  Note that it
was proved in Section~\ref{sec:covering-link} that the numbers
$c_{ik}$ are determined by $r$ and $\sum_k c_{ik}=1$ if $i=2$ or
$\sum_k c_{ik}=0$ otherwise.

If $L$ were concordant to a homology boundary link admitting
pattern~$r$, then the signature jump function of the link $\tilde L =
\bigcup_{i>1,k} i_{c_{ik}} t^k \tilde K_i$ constructed as above would
have the period~$2\pi$, since the signature jump function of $\tilde
L$ is a concordance invariant of $L$.  Let $L$ be the distant union of
$(n-2)$-component unlink and a link whose Seifert matrix and pattern
are as in Theorem~\ref{thm:link-not-concordant-to-boundary-link}. We
will show that for any sufficiently large~$m$, $L$ does not satisfy
the above periodicity condition.  This completes the proof of
Theorem~\ref{thm:link-not-concordant-to-hbl-pattern}.

We assume $m>0$ and fix $p=3$. Let $c_{ik}$ be the numbers determined
by the pattern $r$ and let $\tilde L = \bigcup_{i,k} i_{c_{ik}} t^k
\tilde K_i$ as above. Since the split unlink part has no contribution
to the signature, the signature of $\tilde L$ is equal to the
signature of $\bigcup_{k} i_{c_{2k}} t^k \tilde K_2$. By
Theorem~\ref{thm:covering-seifert-matrix-2comp}, we have $
\delta_{\tilde L}(\theta) = \delta^q_V(y_1\theta)+
\delta^q_V(y_2\theta)+ \delta^q_V(y_0\theta) $ where $y_0 = x_0-x_2$,
$y_1 = x_1-x_0$, $y_2 = x_2-x_1$ and $\{x_i\}$ is a solution of
$$
\begin{bmatrix}
m & 1-m & 0 \\
0 & m & 1-m \\
1-m & 0 & m
\end{bmatrix}
\begin{bmatrix} x_0 \\ x_1 \\ x_2 \end{bmatrix}
= \begin{bmatrix} c_{20} \\ c_{21} \\ c_{22} \end{bmatrix}
$$
as in the proof of
Theorem~\ref{thm:link-not-concordant-to-boundary-link}.  By solving
the equations, we have $y_i = (a_i m + b_i)/(3m^2-3m+1)$, where
\begin{alignat*}{2}
a_0 &= 3c_{20}+3c_{21}-2,  &\quad b_0 &= 1-c_{20}-2c_{21},\\
a_1 &= 1-3c_{20},     &\quad b_1 &= 2c_{20}+c_{21}-1,\\
a_2 &= 1-3c_{21},     &\quad b_2 &= c_{20}-c_{21}.
\end{alignat*}

Since $a_0\equiv a_1\equiv a_2\equiv 1\pmod3$ and $a_0+a_1+a_2=0$, we
may assume that $|a_0|>|a_1|,|a_2|$ by permuting indices.  Choose
minimal $\theta_0>0$ such that $\delta^q_V(\theta_0)\ne 0$. Choose
$\epsilon>0$ such that $\delta^q_V(\theta)=0$ for all
$0<|\theta-\theta_0|<\epsilon$.  Since $\lim_{m\to\infty} |y_i/y_0| <
1$ for $i=1,2$, we can choose $\epsilon'>0$ such that $|y_1/y_0|,
|y_2/y_0| < 1-\epsilon'$ for any large $m$.  We remark that
$0<2\pi|y_0|<\epsilon$ and $2\pi|y_1|, 2\pi|y_2| < \epsilon'\theta_0$
are satisfied for any large $m$ since $y_i \to 0$ as $m\to \infty$.

Let $\theta_1 = \theta_0/y_0$. We claim that for any large $m$,
$\delta_{\tilde L}(\theta_1) \ne 0$ and $\delta_{\tilde
L}(\theta_1+2\pi)=0$. Since $|y_1/y_0|, |y_2/y_0| < 1$,
$\delta_{\tilde L}(\theta_1) = \delta^q_V(\theta_0) \ne 0$. Since
$0<|2\pi y_0|<\epsilon$ and $|\theta_0 y_i/y_0+2\pi y_i| < \theta_0$
for $i=1,2$, $\delta_{\tilde L}(\theta_1+2\pi) =
\delta^q_V(\theta_0+2\pi y_0) + \delta^q_V(\theta_0 y_1/y_0+2\pi y_1)
+ \delta^q_V(\theta_0 y_2/y_0+2\pi y_2) = 0$. This proves the claim.

The claim implies that $\delta_{\tilde L}(\theta)$ is not of
period~$2\pi$.  Therefore $L$ is not concordant to any homology
boundary links admitting pattern~$r$ if $m$ is sufficiently large.
\end{proof}

\end{document}